\documentclass[11pt]{article}
\usepackage{amsmath,amssymb,amsthm,url}
\usepackage[letterpaper,hmargin=1.25in,vmargin=1in]{geometry}
\usepackage{tikz}
\usetikzlibrary{matrix}

\newtheorem{thm}{Theorem}

\newtheorem{prop}[thm]{Proposition}
\newtheorem{cor}[thm]{Corollary}

\theoremstyle{definition}
\newtheorem{defn}[thm]{Definition}
\newtheorem{example}[thm]{Example}

\newcommand{\coP}{\overline{P_5}}

\title{Induced paths in strongly regular graphs}
\author{Robert F.\ Bailey\footnote{Corresponding author.  School of Science and the Environment (Mathematics), Memorial University--Grenfell Campus, 20 University Drive, Corner Brook, NL A2H~6P9, Canada. E-mail: \texttt{rbailey@grenfell.mun.ca}.}\,\, and Abigail K.\ Rowsell\footnote{School of Science and the Environment (Mathematics), Memorial University--Grenfell Campus, 20 University Drive, Corner Brook, NL A2H~6P9, Canada. E-mail: \texttt{abigailrowsell21@gmail.com}.}}

\begin{document}

\maketitle

\begin{abstract}
This paper studies induced paths in strongly regular graphs.  We give an elementary proof that a strongly regular graph contains a path $P_4$ as an induced subgraph if and only if it is primitive, i.e.\ it is neither a complete multipartite graph nor its complement.  Also, we investigate when a strongly regular graph has an induced subgraph isomorphic to $P_5$ or its complement, considering several well-known families including Johnson and Kneser graphs, Hamming graphs, Latin square graphs, and block-intersection graphs of Steiner triple systems.\\


\end{abstract}

\section{Introduction}
In this paper, all graphs are finite, with no loops or multiple edges.  In particular, we are concerned with the following class of graphs.

\begin{defn} \label{defn:SRG}
A graph $G$ is {\em strongly regular} with parameters $(n,k,\lambda,\mu)$ if it has $n$ vertices, is regular with degree $k$, any two adjacent vertices have $\lambda$ common neighbours, and any two non-adjacent vertices have $\mu$ common neighbours.
\end{defn}

It is well-known that $G$ is strongly regular with parameters $(n,k,\lambda,\mu)$ if and only if its complement $\overline{G}$ is strongly regular with parameters $(n,\overline{k},\overline{\lambda},\overline{\mu})$, where $\overline{k}=n-k-1$, $\overline{\lambda}=n-2k+\mu-2$ and $\overline{\mu}=n-2k+\lambda$.  From the definition, it follows that a strongly regular graph has diameter~$2$ unless $\mu=0$, whereby it must be disconnected.

A strongly regular graph is called {\em primitive} if it and its complement is both connected, and {\em imprimitive} otherwise; it is well-known that the only imprimitive, connected strongly regular graphs are the complete multipartite graphs $K_{r\times m}$ (with $r$ parts of size $m$), which have parameters $(rm,\ (r-1)m,\ (r-2)m, (r-1)m)$; their complements are the disjoint union of $r$ copies of $K_m$.  More information on strongly regular graphs can be found in~\cite{BrouwerHaemers12,BvM}.

In this paper, we are concerned with strongly regular graphs with, our without, specific induced subgraphs.  The following definition is relevant here.

\begin{defn} \label{defn:Hfree}
A graph $G$ is called {\em $H$-free} if it contains no induced subgraph isomorphic to some graph $H$, and is called {\em $\mathcal{F}$-free} if it contains no induced subgraph isomorphic to a member of some family of graphs $\mathcal{F}$.
\end{defn}

A particularly important class of such graphs is as follows.

\begin{defn} \label{defn:cograph}
A graph $G$ is a {\em cograph} if it contains no subgraph isomorphic to a path on $4$ vertices, $P_4$.
\end{defn}

The name ``cograph'' is a shorthand for ``complement reducible graph'', because of the equivalence between graphs which are $P_4$-free and graphs which may be reduced to isolated vertices by recursively complementing all connected subgraphs.  This equivalence forms part of the ``Fundamental Theorem on Cographs'' of Corneil, Lerchs and Stewart Burlingham~\cite[Theorem~2]{Corneil81} (see also~\cite[Theorem~11.3.3]{classes}).  Nowadays, having no induced $P_4$ is frequently taken as the definition of the term ``cograph'' (see, for instance, \cite{Cameron21,Epple15,Johnson15}).

Clearly, a connected $P_4$-free graph necessarily has diameter at most~$2$, as is the case with strongly regular graphs, so we begin with an examination of how these classes coincide.

\section{Induced $4$-paths} \label{sec:P4}
It is a simple exercise to see that a complete multipartite graph has no induced subgraphs isomorphic to $P_4$; thus, any imprimitive strongly regular graph is a cograph.  Conversely, it follows from the notion of complement-reducibility that a connected cograph must have a disconnected complement, which shows that a primitive strongly regular graph is not a cograph and thus contains $P_4$ as an induced subgraph.  The purpose of this section is to give an elementary proof of this fact, using only the properties of strongly regular graphs, which is the following result.

\begin{thm} \label{thm:P4}
Let $G$ be a primitive strongly regular graph with parameters $(n,k,\lambda,\mu)$.  Then $G$ contains an induced subgraph isomorphic to $P_4$.
\end{thm}

\proof We consider different values for the parameters $\lambda$ and $\mu$.  Since we are considering primitive strongly regular graphs only, we know that $0<\mu<k$.  For any vertex $u$ of a graph $G$, we let $G_i(u)$ denote the subset of vertices at distance $i$ from $u$.

The easiest case is when $\lambda=0$ and $\mu=1$, the so-called Moore graphs.  In this case, $G$ has girth $5$, and so clearly contains a $5$-cycle as an induced subgraph, and thus also an induced $P_4$.

Next, we suppose that $\lambda=0$ and $\mu>1$; in this case, $G$ has girth $4$.  Choose vertices $u,v,w$, where $v\in G_1(u)$ and $w\in G_2(u)\cap G_1(v)$, so that $uvw$ is an induced $P_3$.  Now choose some $x\in G_1(u)\setminus G_1(w)$: since $\mu<k$ we know that $G_1(u)\setminus G_1(w)$, which has size $k-\mu$, is non-empty.  Then $xuvw$ is an induced $P_4$ in $G$.

It remains to consider the graphs with girth $3$, i.e.\ those with $\lambda>0$.  First, we suppose that $\mu\leq\lambda+1$.  Again, we choose vertices $u,v,w$ such that $v\in G_1(u)$ and $w\in G_2(u)\cap G_1(v)$.  The aim is to construct an induced $P_4$ of the form $uvwx$; the new vertex $x$ must be a neighbour of $w$, must lie in $G_2(u)$ (as otherwise, $uvwx$ would contain a $4$-cycle), and must not be a neighbour of $v$ (as otherwise, $uvwx$ would contain a $3$-cycle).  Now, $w$ has $k$ neighbours, of which $\mu$ lie in $G_1(u)$ and $k-\mu$ lie in $G_2(u)$.  Also, $v$ has $k$ neighbours: as well as $u$ and $w$, these are $\lambda$ vertices in $G_1(u)$, and a further $k-\lambda-2$ vertices in $G_2(u)$.  Since $\mu\leq\lambda+1$, we have $k-\mu>k-\lambda-2$, so by the pigeonhole principle, the set $(G_2(u)\cap G_1(w)) \setminus (G_2(u)\cap G_1(v))$ is non-empty.  Hence, a suitable vertex $x$ exists so that $uvwx$ is an induced $P_4$ in $G$.

Finally, we suppose that $\mu\geq\lambda+1$.  This time, let $u,v,w$ be vertices such that $w\in G_1(u)$, and $v$ is not adjacent to either $u$ or $w$.  Our aim is to find a vertex $x\in G_1(u)\cap G_1(v)$ but where $x\not\in G_1(w)$, so that $wuxv$ is an induced $P_4$ in $G$.  Now, $|G_1(u)\cap G_1(v)|=\mu$, since this is precisely the set of common neighbours of the non-adjacent vertices $u,v$.  Also, $w$ has $\lambda$ neighbours in $G_1(u)$, so $x$ cannot be any of these vertices.  However, since $\mu>\lambda$, there exists a suitable vertex $x$ so that $wuxv$ is an induced $P_4$, as required.

This concludes the proof.  \endproof

We remark that the latter part of the proof (where $\mu\geq\lambda+1$) is really just applying the previous argument to the complementary graph $\overline{G}$, and constructing an induced $P_4$, namely $uvwx$, in $\overline{G}$ in exactly the same manner as is done in that part of the proof; since $P_4$ is self-complementary, the complement of the induced $P_4$ in $\overline{G}$ is an induced $P_4$ in $G$.  Also, either method can be used when $\mu=\lambda+1$, which is reassuring, given that if $G$ and $\overline{G}$ have the same parameters we necessarily have $\mu=\lambda+1$ in that case.

As a consequence of the Theorem above, we have the following fact.

\begin{cor} \label{cor:SRcograph}
A strongly regular graph $G$ is $P_4$-free if and only if it is imprimitive (i.e.\ a complete multipartite graph or its complement).
\end{cor}

\section{Induced $5$-paths} \label{sec:P5}

Given the straightforward characterization of $P_4$-free strongly regular graphs, it seems natural to consider extensions of this question.  Perhaps the most natural next step is to consider those graphs which are $P_5$-free, or those which are $\{P_5,\overline{P_5}\}$-free.  For further details on this latter class, see Chudnovsky {\em et al.}~\cite{Chudnovsky17}, where an algorithmic characterization of such graphs is obtained, analogous to that obtained by Corneil, Perl and Stewart for cographs~\cite{Corneil85}.

Clearly, if a graph is $P_4$-free, then it cannot contain an induced $P_5$ or $\overline{P_5}$ either.  However, there are examples of primitive strongly regular graphs which are $P_5$-free or $\overline{P_5}$-free: for example, the Petersen graph has girth~$5$, so has no induced $\overline{P_5}$ (which is formed of a $3$-cycle and a $4$-cycle with a edge in common); consequently, its complement (the Johnson graph $J(5,2)$) contains no induced $P_5$.  However, as we will see in the subsections below, many well-known families of primitive strongly regular graphs do contain both induced $P_5$ and $\overline{P_5}$ subgraphs.

\subsection{Johnson, Kneser and Hamming graphs}
The {\em Johnson graph} $J(m,2)$, also known as the {\em triangular graph} $T(m)$, has as its vertices the $2$-subsets of $\{1,\ldots,m\}$, and two $2$-subsets are adjacent if their intersection has size~$1$; this graph is strongly regular with parameters $\left( \binom{m}{2}, 2(m-2), m-2, 4 \right)$.  The complement of $J(m,2)$ is the {\em Kneser graph} $K(m,2)$.  Also, the {\em Hamming graph} $H(2,m)$, also known as the {\em square lattice graph}, has all ordered pairs of symbols from $\{0,\ldots,m-1\}$ as its vertices, with two pairs adjacent whenever they agree in a single co-ordinate; it is strongly regular with parameters $(m^2, 2(m-1), m-2, 2)$.

\begin{prop} \label{prop:JK}
The Johnson graph $J(m,2)$ contains an induced $P_5$ if and only if $m\geq 6$, and an induced $\coP$ if and only if $m\geq 5$.
\end{prop}

\proof For $m\leq 3$, $J(m,2)$ has fewer than $5$ vertices, while $J(4,2)$ is the complete multipartite graph $K_{2,2,2}$.  So for $m\leq 4$, $J(m,2)$ is $\{P_5,\coP\}$-free.  For $m\geq 6$, it is straightforward to verify that the following is an induced $P_5$ subgraph of $J(m,2)$,
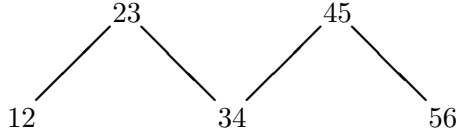
\begin{center}
\setlength{\unitlength}{7mm}
\begin{picture}(10,4)
\thicklines
\put(0.75,0.8){12}
\put(2.75,2.8){23}
\put(4.75,0.8){34}
\put(6.75,2.8){45}
\put(8.75,0.8){56}

\put(1.3,1.3){\line(1,1){1.4}}
\put(3.3,2.7){\line(1,-1){1.4}}
\put(5.3,1.3){\line(1,1){1.4}}
\put(7.3,2.7){\line(1,-1){1.4}}
\end{picture}
\end{center}
while for $m\geq 5$, we have the following induced $P_5$ in $K(m,2)$ (which yields an induced $\coP$ in $J(m,2)$).
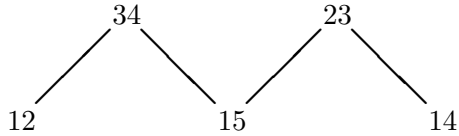
\begin{center}
\setlength{\unitlength}{7mm}
\begin{picture}(10,4)
\thicklines
\put(0.75,0.8){12}
\put(2.75,2.8){34}
\put(4.75,0.8){15}
\put(6.75,2.8){23}
\put(8.75,0.8){14}

\put(1.3,1.3){\line(1,1){1.4}}
\put(3.3,2.7){\line(1,-1){1.4}}
\put(5.3,1.3){\line(1,1){1.4}}
\put(7.3,2.7){\line(1,-1){1.4}}
\end{picture}
\end{center}
Finally, the Petersen graph $K(5,2)$ cannot contain an induced $\coP$ (as we observed earlier), and thus $J(5,2)$ has no induced $P_5$.  This completes the proof. \endproof

\begin{prop} \label{prop:Hamming}
The Hamming graph $H(2,m)$ contains an induced $P_5$ and an induced $\coP$ if and only if $m\geq 3$.
\end{prop}

\proof For $m\leq 2$, $H(2,m)$ has fewer than five vertices.  For $m\geq 3$, we have the following induced $P_5$ and $\coP$ subgraphs.
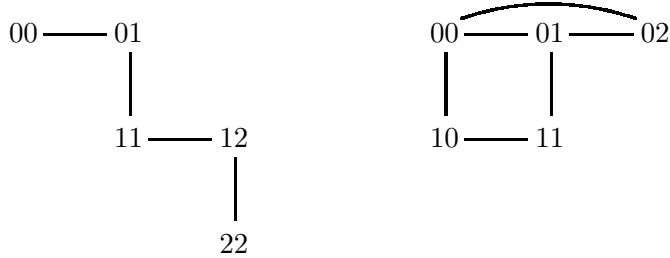
\begin{center}
\setlength{\unitlength}{7mm}
\begin{picture}(14,6)
\thicklines
\put(0.75,4.8){00}
\put(2.75,4.8){01}
\put(2.75,2.8){11}
\put(4.75,2.8){12}
\put(4.75,0.8){22}

\put(1.4,4.95){\line(1,0){1.2}}
\put(3.4,2.95){\line(1,0){1.2}}
\put(3.05,3.4){\line(0,1){1.2}}
\put(5.05,1.4){\line(0,1){1.2}}


\put(8.75,4.8){00}
\put(10.75,4.8){01}
\put(12.75,4.8){02}
\put(8.75,2.8){10}
\put(10.75,2.8){11}

\put(9.4,4.95){\line(1,0){1.2}}
\put(11.4,4.95){\line(1,0){1.2}}
\put(9.4,2.95){\line(1,0){1.2}}
\put(9.05,3.4){\line(0,1){1.2}}
\put(11.05,3.4){\line(0,1){1.2}}
\qbezier(9.35,5.25)(11,5.8)(12.65,5.25)

\end{picture}
\end{center}
\vspace{-7mm}
\endproof

\subsection{Latin square graphs}
Recall that a {\em Latin square} of order $m$ is an $m\times m$ array filled with symbols from the set $\{0,\ldots,m-1\}$, so that each symbol occurs exactly once in each row and once in each column.  From a Latin square, the {\em Latin square graph} has $m^2$ vertices corresponding to the cells of the array, with two cells being adjacent whenever they are in the same row, are in same column, or are filled with the same symbol.  These are well-known to be strongly regular with parameters $(m^2,3(m-1),m,6)$.  The next two results give sufficient conditions for such a graph to have induced $P_5$ and $\coP$ subgraphs.

\begin{prop} \label{prop:LS}
Let $L$ be a Latin square of order $m\geq 5$.  The then corresponding Latin square graph contains an induced $P_5$ subgraph.
\end{prop}

\proof We will directly construct a suitable $P_5$ subgraph.  Without loss of generality, we can assume that the first row of $L$ contains the symbols $0,1,\ldots,m-1$ in order, so we choose the first two cells as the first two vertices of our path.  Next, the second column must contain the symbol $2$ in some cell, so we choose this cell as the third vertex (by permuting rows if necessary, we may assume that this cell is in the second row).  The next vertex will also be in the second row, but it cannot be in the first column, nor can it contain the symbols $0$ or $1$; of the remaining $m-1$ cells in that row, at most three are unavailable, but since $m\geq 5$ we have that $m-4>0$, so at least one cell will be available.  Choose one of these cells, and suppose that symbol $3$ occurs in it.  For the fifth vertex, we will choose a cell which also contains symbol $3$.  We cannot choose a cell from the first row or the first two columns, which leaves $m-4$ possible cells to choose from, and since $m-4>0$ we are done.
\endproof

\begin{prop} \label{prop:LS-coP5}
Let $L$ be a Latin square of order $m\geq 6$.  The then corresponding Latin square graph contains an induced $\coP$ subgraph.
\end{prop}

\proof As with Proposition~\ref{prop:LS}, we will construct a suitable $\coP$ directly.  Again, we assume that the first row of $L$ contains the symbols $0,1,\ldots,m-1$ in order; we choose the first three cells (with symbols $0,1,2$) to form the $3$-cycle in our $\coP$.  For the remaining two vertices, we will choose two which are in the first two columns and share a row.  However, we cannot choose cells containing symbols $1$ or $2$ from the first column, or $0$ or $2$ from the second column, which means that, of the $m-1$ remaining rows, at most four are unavailable.  But since $m\geq 6$, we have $m-5>0$, so there must be at least one row containing two new symbols in the first two columns.  We choose the first two cells in such a row, and we are done. \endproof

We demonstrate these methods in the following example.

\begin{example} \label{example:LS6}
In the cyclic Latin square of order $6$ as shown below, the highlighted cells on the left yield an induced $P_5$, while those on the right yield an induced $\coP$.
\begin{center}
\begin{tikzpicture}

\tikzset{square matrix/.style={
    matrix of nodes,
    column sep=-\pgflinewidth, row sep=-\pgflinewidth,
    nodes={draw,
      minimum height=#1,
      anchor=center,
      text width=#1,
      align=center,
      inner sep=0pt
    },
  },
  square matrix/.default=0.7cm
}

\hspace{-3cm}
\matrix[square matrix]
{
|[fill=lightgray]|0 & |[fill=lightgray]|1 & 2 & 3 & 4 & 5 \\
1 & |[fill=lightgray]|2 & |[fill=lightgray]|3 & 4 & 5 & 0 \\
2 & 3 & 4 & 5 & 0 & 1 \\
3 & 4 & 5 & 0 & 1 & 2 \\
4 & 5 & 0 & 1 & 2 & 3 \\
5 & 0 & 1 & 2 & |[fill=lightgray]|3 & 4 \\
};

\hspace{6cm}
\matrix[square matrix]
{
|[fill=lightgray]|0 & |[fill=lightgray]|1 & |[fill=lightgray]|2 & 3 & 4 & 5 \\
1 & 2 & 3 & 4 & 5 & 0 \\
2 & 3 & 4 & 5 & 0 & 1 \\
|[fill=lightgray]|3 & |[fill=lightgray]|4 & 5 & 0 & 1 & 2 \\
4 & 5 & 0 & 1 & 2 & 3 \\
5 & 0 & 1 & 2 & 3 & 4 \\
};

\end{tikzpicture}
\end{center}
\end{example}

What happens for the remaining values of $m$?  For $m\leq 2$, there are fewer than five vertices, while for $m=3$ the graph which arises is a complete multipartite graph, so no induced $P_5$ or $\coP$ is possible.  For $m=4$ and $m=5$, there are exactly two {\em main classes} (see~\cite[{\S}III.1]{handbook}) of Latin squares of each order, and thus two non-isomorphic Latin square graphs for these orders; a slight modification of the method is required, but induced $P_5$ and $\coP$ subgraphs can be found in each case, which is left as an exercise for the reader, using the examples in~\cite[{\S}III.1]{handbook}.

For the next family, we recall that two Latin squares $L,M$ of the same order are {\em orthogonal} if, when superimposed, each ordered pair of symbols occurs in exactly one cell of the array.  From this, we obtain a strongly regular graph (also known as a Latin square graph) on the set of cells, where two cells are adjacent if they are in the same row, in the same column, share a symbol in $L$, or share a symbol in $M$; the parameters are $\left(m^2, 4(m-1), m+4, 12\right)$.\footnote{In general, given a collection of $t$ mutually orthogonal Latin squares, or MOLS, of order $m$, a {\em Latin square graph} has the $m^2$ cells as vertices, and two cells are adjacent if they are in the same row or column, or share a symbol in one of the Latin squares. This graph has parameters $\left(m^2, (t+2)(m-1), m-2+t(t+1), (t+1)(t+2) \right)$; a graph with these parameters, but not necessarily arising from MOLS, is a {\em pseudo-Latin square graph}.}

\begin{prop} \label{prop:MOLS}
Let $L$ and $M$ be orthogonal Latin squares of order $m\geq 8$.  The then corresponding Latin square graph contains an induced $P_5$ subgraph.
\end{prop}

\proof We will construct an induced $P_5$ in a similar manner to the proof of Proposition~\ref{prop:LS}.  Without loss of generality, we can assume that the first rows of $L$ and $M$ both contain the symbols $0,1,\ldots,m-1$ in order.  Thinking of the vertices as the cells of an $m\times m$ array whose entries are ordered pairs of symbols, we assume that the first two vertices of our path are the cells $00$ and $11$ in the first row, and that the third vertex is in the second row and contains two new symbols, e.g.\ $23$.  For the fourth vertex, in the second row we take a vertex which cannot have $0$ or $1$ in either coordinate, and cannot be in the first two columns, so there are at most six cells we must avoid.  But since $m\geq 8$, we have $m-6>0$, so there is a cell we can choose.

Suppose without loss of generality that the fourth vertex is $45$.  For the fifth vertex, we will choose a vertex with $4$ in the first co-ordinate, that is not adjacent to the first three vertices.  Taking rows, columns and entries into consideration, there are at most seven cells we must avoid, but since $m\geq 8$, a suitable cell is guaranteed to exist. \endproof

\begin{prop} \label{prop:MOLS-coP5}
Let $L$ and $M$ be orthogonal Latin squares of order $m\geq 10$.  The then corresponding Latin square graph contains an induced $\coP$ subgraph.
\end{prop}

\proof The proof is very similar to that of Proposition~\ref{prop:LS-coP5}: we label the vertices in the same manner as in Proposition~\ref{prop:MOLS}, choose vertices $00$, $11$ and $22$ from the first row, and two vertices from the first two columns in the same row.  This time, there are at most eight other rows we must avoid, but since $m\geq 10$, a suitable row exists, and thus we are able to construct an induced $\coP$.  \endproof

We remark that the ``threshold'' values for $m$ which appear in the proofs of Propositions~\ref{prop:LS} to~\ref{prop:MOLS-coP5} are probably artificially high; it may be possible to adapt the proofs to cover more cases.  Also, the method of proof appears as if it should generalize, i.e.\ that for each $t$, there should be some value $N_t(m)$ such that for all $m\geq N_t(m)$, any Latin square graph arising from a set of $t$ mutually orthogonal Latin squares of order $m$ contains an induced $P_5$, and likewise for induced $\coP$ subgraphs.


\subsection{Block-intersection graphs of Steiner triple systems}
A {\em Steiner triple system} of order $m$, or $\mathrm{STS}(m)$, is a pair $(X,\mathcal{B})$ where $X$ is a set of $m$ points, and $\mathcal{B}$ is a collection of $3$-subsets of $X$, called blocks, with the property that any pair of points from $X$ lies in exactly one block in $\mathcal{B}$.  It is well-known that a Steiner triple system exists if and only if $m\equiv 1$ or $3$~mod~$6$.  There are unique Steiner triple systems of orders $3$, $7$ and $9$, two of order $13$, and $80$ of order $15$ (see~\cite{Mathon83}); for larger orders the number of systems is in the billions.  For more information, see~\cite[{\S}II.2]{handbook}.

The {\em block-intersection graph} of a Steiner triple system has vertex set $\mathcal{B}$, and two blocks are adjacent if their intersection is non-empty.  It is straightforward to show that this graph is strongly regular, with parameters $\left( \frac{1}{6}m(m-1),\, \frac{3}{2}(v-3), \, \frac{1}{2}(v+3),\, 9 \right)$.  We note that the block-intersection graphs of the unique ${\rm STS}(3)$, ${\rm STS}(7)$ and ${\rm STS}(9)$ are $K_1$, $K_7$ and $K_{3,3,3,3}$ respectively, none of which can contain an induced $P_5$ or $\overline{P_5}$.  

For the block-intersection graphs of the two $\mathrm{STS}(13)$s, we have induced $P_5$ and $\coP$ subgraphs as follows.

\begin{example} \label{example:STS13}
The blocks of the two $\mathrm{STS}(13)$s, with point set $\{1,\ldots 13\}$, are given below:
\[ \begin{array}{ccccccc} 
1\, 2\, 3 \,  & \, 1\, 4\, 5 \,  & \, 1\, 6\, 11 \, & \, 1\, 7\, 8 \,   & \, 1\, 9\, 10 \,  & \, 1\, 12\, 13 \, & \, 2\, 4\, 8 \\
2\, 5\, 7 \,  & \, 2\, 6\, 10 \, & \, 2\, 9\, 12 \, & \, 2\, 11\, 13 \, & \, 3\, 4\, 11 \,  & \, 3\, 5\, 10 \,  & \, 3\, 6\, 12 \\
3\, 7\, 9 \,  & \, 3\, 8\, 13 \, & \, 4\, 6\, 7 \,  & \, 4\, 9\, 13 \,  & \, 4\, 10\, 12 \, & \, 5\, 6\, 13 \,  & \, 5\, 8\, 12 \\
5\, 9\, 11 \, & \, 6\, 8\, 9 \,  & \, 7\, 10\, 13\, & \, 7\, 11\, 12 \, & \, 8\, 10\, 11 
\end{array} \]
and
\[ \begin{array}{ccccccc} 
1\, 2\, 3 \,  & \, 1\, 4\, 5 \,  & \, 1\, 6\, 11 \, & \, 1\, 7\, 8 \,   & \, 1\, 9\, 10 \,  & \, 1\, 12\, 13 \, & \, 2\, 4\, 8 \\
2\, 5\, 9 \,  & \, 2\, 6\, 10 \, & \, 2\, 7\, 13 \, & \, 2\, 11\, 12 \, & \, 3\, 4\, 11 \,  & \, 3\, 5\, 10 \,  & \, 3\, 6\, 12 \\
3\, 7\, 9 \,  & \, 3\, 8\, 13 \, & \, 4\, 6\, 7 \,  & \, 4\, 9\, 12 \,  & \, 4\, 10\, 13 \, & \, 5\, 6\, 13 \,  & \, 5\, 7\, 11 \\
5\, 8\, 12 \, & \, 6\, 8\, 9 \,  & \, 7\, 10\, 12\, & \, 8\, 10\, 11 \, & \, 9\, 11\, 13
\end{array} \]
Both of these have $1\, 2\, 3$,\, $1\, 4\, 5$,\, $4\, 6\, 7$,\, $6\, 8\, 9$,\, $8\, 10\, 11$ as an induced $P_5$ in their block-intersection graphs.  Also, they each have $1\, 2\, 3$,\, $1\, 4\, 5$,\, $1\, 9\, 10$,\, $2\, 4\, 8$,\, $4\, 6\, 7$ as an induced $\coP$.
\end{example}

For the $80$ distinct $\mathrm{STS}(15)$s, we were able to verify the existence of induced $P_5$ subgraphs computationally, using the GAP computer algebra system~\cite{GAP}: first, we constructed the 80 systems using the DESIGN package~\cite{DESIGN}, and obtained their block-intersection graphs with the GRAPE package~\cite{GRAPE}.  Then for each graph, by enumerating $5$-subsets of vertices and examining the corresponding induced subgraphs, we could quickly verify (in GRAPE) the existence of an induced $P_5$ in each of them.  (We will see below that this computation was unnecessary for induced $\coP$ subgraphs.)

In our next two results, we will use the following terminology: in an ${\rm STS}(m)$, given a pair of distinct points $x,y$, we call the unique block containing $x$ and $y$ the {\em completion} of $x\, y$.

\begin{thm} \label{thm:STS-P5}
The block intersection graph of a Steiner triple system of order $m$ contains an induced $P_5$ if and only if $m\geq 13$.
\end{thm}

\proof We already known that the block-intersection graph of an ${\rm STS}(m)$ is $P_5$-free for $m<13$, and contains an induced $P_5$ for $m=13$ and $m=15$; from now on, we will assume that $m\geq 19$.

Let $\mathcal{S}$ be an ${\rm STS}(m)$, with point set $\{1,\ldots,m\}$ and block set $\mathcal{B}$.  Without loss of generality, we can assume that $\mathcal{S}$ contains blocks $A=1\,2\,3$ and $B=1\,4\,5$.  We can also assume that there is a block $C=4\,6\,7$: there are $\frac{m-3}{2}$ blocks of the form $4\,i\,j$ distinct from $B$ (completions of the pair $4\,i$ where $i\not\in\{1,2,3\}$), but only two of these can intersect with $A$; since $m\geq 19$ we must have $\frac{m-3}{2}>2$, so a suitable block exists.

We will find blocks $D$ and $E$ which satisfy the following: (i) $D$ intersects $C$ but neither $A$ nor $B$, and (ii) $E$ intersects $D$ but none of $A$, $B$ or $C$.  (These will yield an induced $P_5$ in the block-intersection graph of $\mathcal{S}$.)  Without loss of generality, assume that $6\in D$.  Now, there are $\frac{m-3}{2}$ blocks containing $6$ other than $C$; these include the completions of $1\,6$, $2\,6$, $3\,6$ and $5\,6$, which intersect with $A$ or $B$.  So there are at most four blocks which we cannot choose as $D$; however, since $m\geq 19$ we have that $\frac{m-3}{2}>4$, so such a block $D$ must exist.  By relabelling points if necessary, we may assume that $D=6\,8\,9$.

The argument to find $E$ is similar: assume that $8\in E$.  There are $\frac{m-3}{2}$ blocks containing $8$ other than $D$; this time, to avoid intersecting with $A$, $B$ or $C$, we cannot use the completions of $1\,8, \ldots, 5\,8$ or $7\,8$, which yield at most six blocks.  However, since $m\geq 19$ we have that $\frac{m-3}{2}>6$, and such a block exists.  Again by relabelling, we may assume that $E=8\,10\,11$. \endproof

The proof for the existence of induced $\coP$ subgraphs is similar, but this time we do not need to treat $m=13$ or $15$ separately.

\begin{thm} \label{thm:STS-coP5}
The block intersection graph of a Steiner triple system of order $m$ contains an induced $\overline{P_5}$ if and only if $m\geq 13$.
\end{thm}

\proof Let $\mathcal{S}$ be an ${\rm STS}(m)$ as in the proof of Theorem~\ref{thm:STS-P5}; again we can assume that there are blocks $A=1\,2\,3$, $B=1\,4\,5$ and $C=4\,6\,7$ in $\mathcal{S}$.  This time, we will obtain blocks $D$ and $E$ which satisfy the following: (i) $D$ intersects $A$ and $C$ but not $B$, and (ii) $E$ intersects $A$ and $B$ but neither $C$ nor $D$.  (These will yield an induced $\overline{P_5}$ in the block-intersection graph of $\mathcal{S}$.)

For a block $D$ to satisfy (i), it must be the completion of a pair $2\,6$, $2\,7$, $3\,6$ or $3\,7$; i.e.\ a block of the form $2\,6\,x$, $2\,7\,y$, $3\,6\,z$ or $3\,7\,w$ for some points $x,y,z,w$.  Now, we cannot have $x,y,z,w\in\{1,2,3,4,6,7\}$, as this would repeat a pair that appears in $A$, $B$ or $C$.  Also, we cannot have $x=y=z=w=a_5$, as this would cause some of the pairs $2\,5$, $3\,5$, $5\,6$ or $5\,7$ to be repeated.  Consequently, at least one of these completions must use a new point; without loss of generality we may assume that $x=8$, i.e.\ that $D=2\,6\,8$.

To obtain a block $E$ which satisfies (ii), we show that there exists a block $1\,a\,b$ with the desired property.  There are $\frac{m-1}{2}$ blocks containing the point $1$, namely $A$, $B$ and $\frac{m-5}{2}$ others.  These others include the completions of $1\,6$, $1\,7$ or $1\,8$, which we cannot choose as $E$.  But since $m\geq 13$, we must have that $\frac{m-5}{2}>3$, and so a suitable block must exist.  By relabelling points if necessary, we may assume that $E=1\,9\,10$. \endproof

\section{Conclusion}
It would be desirable to obtain a complete characterization of the primitive strongly regular graphs which contain an induced $P_5$ or $\coP$, analogous to Theorem~\ref{thm:P4} for induced $P_4$ subgraphs.  As we have already seen, if a graph contains no triangles, it cannot have an induced $\coP$; for strongly regular graphs, such a graph has $\lambda=0$, and there are exactly seven examples known: the $5$-cycle $C_5$, and the Petersen, Clebsch, Hoffman--Singleton, Gewirtz, $M_{22}$ and Higman--Sims graphs, which have $5$, $10$, $16$, $50$, $56$, $77$ and $100$ vertices respectively.  While none of these can contain an induced $\coP$, using GAP it can be shown that (other than $C_5$) all have an induced $P_5$.

Also, using GAP (and the libraries available at \url{www.distanceregular.org}), we were able to test all primitive strongly regular graphs on up to $28$ vertices; apart from $C_5$, all have an induced $P_5$, and the only examples without an induced $\coP$ are the triangle-free examples mentioned above.  So perhaps it is the case that $\lambda>0$ is sufficient for a strongly regular graph to have an induced $\coP$, and having $\overline{\lambda}>0$ (i.e.\ $\lambda$ in the complement graph) is sufficient to have an induced $P_5$?

Of course, there are plenty of other induced subgraphs one could look for (or look to exclude), such as longer paths, or the ``gem'' on five vertices formed by taking a $P_4$ and adding a new vertex adjacent to all others.  (The classes of $\{\textnormal{gem},\ \overline{\textnormal{gem}}\}$-free and $\{P_5,\textnormal{gem}\}$-free graphs have been of recent interest, for example in~\cite{Chudnovsky20,Karthick18}.)


\subsection*{Acknowledgements}
The first author would like to thank Ortrud Oellermann for introducing him to cographs.  He is supported by an NSERC Discovery Grant.  The second author was supported by an NSERC Undergraduate Student Research Award.


\end{document}